\documentclass{siamltex}
\usepackage{amssymb,amsmath,amscd,verbatim}
\usepackage{enumerate}
\usepackage{algorithmic}
\usepackage{algorithm}
\usepackage{subfigure}

\usepackage{color}
\usepackage[pdftex]{graphicx}
\usepackage{cite}
\newcommand{\bi}{\begin{itemize}}
\newcommand{\ei}{\end{itemize}}

\newcommand{\caliber}{{c}}

\newcommand{\finevar}{F}
\newcommand{\coarsevar}{C}

\newcommand{\Hone}{{\bf H1}}
\newcommand{\Htwo}{{\bf H2}}

\newtheorem{remark}[theorem]{Remark}

\renewcommand{\refname}{References cited}

\title{An algebraic distances measure of AMG strength of connection}
\author{A.~Brandt\footnotemark[2], J.~Brannick\footnotemark[3], K.~Kahl\footnotemark[4], and I.~Livshits\footnotemark[5]}

\date{}
\begin{document}

\maketitle
\begin{abstract}
Algebraic multigrid is an iterative method that is often optimal
for solving the matrix equations that arise in a wide variety of
applications, including discretized partial differential equations. It
automatically constructs a sequence of increasingly smaller matrix
problems that enable efficient resolution of all scales
present in the solution.  One of the main components of the method is an adequate choice of coarse grids. The current coarsening methodology is based on measuring how a
so-called {\em algebraically smooth} error value at one point depends on the error  values at its neighbors.  Such a concept of {\em strength of
connection} is well understood for operators whose principal part is
an M-matrix; however, the strength concept for more general matrices
is not yet clearly understood, and this lack of knowledge limits the
scope of AMG applicability. The purpose of this paper is to motivate a
general definition of strength of connection, based on the notion of {\em algebraic 
distances}, discuss its
implementation, and present the results of initial numerical
experiments.  
The algebraic distance measure, we propose, uses as its main 
tool a least squares functional, which   is also applied to define interpolation.  
 \end{abstract}
 

\begin{AM}
65F10, 65N55, 65F30
\end{AM}
 \footnotetext[2]{Department of Mathematics, University of California Los Angeles, Los Angeles CA, 90095, USA~(abrandt@math.ucla.edu).}
\footnotetext[3]{Department of Mathematics,
Pennsylvania State University, University Park, PA 16802, USA (brannick@psu.edu). \thanks{Brannick's work was supported by the National Science Foundation under grants OCI-0749202 and DMS-810982.}}
\footnotetext[4]{Fachbereich Mathematik und Naturwissenschaften,
Bergische Universit\"at Wuppertal, D-42097 Wuppertal, Germany, 
(kkahl@math.uni-wuppertal.de). \thanks{Kahl's work was supported by the
  Deutsche Forschungsgemeinschaft through the Collaborative Research
  Centre SFB-TR 55 ``Hadron physics from Lattice QCD''}}
\footnotetext[5]{Department of Mathematical Sciences, Ball State University, Muncie IN, 47306, USA (ilivshits@bsu.edu), \thanks{Livshits' work was supported by DOE subcontract B591416 and NSF DMS-0521097.}}
\renewcommand{\thefootnote}{\arabic{footnote}}

\section{Introduction}
Multigrid methods for solving  general systems of linear equations  
$Ax=b$ are all based on the smoothing property of relaxation. For
any $0 < \tau < 1$, an error vector $e$ is called $\tau$-smooth if
all its 
residuals
are smaller than $\tau
\|e\|$. The basic observation \cite{B83} is that the convergence
of a proper relaxation process
 slows down only when the current error is
$\tau$-smooth with $\tau \ll 1$, the smaller the $\tau$ the slower
the convergence. This observation and the assumption that when
relaxation slows down, the error vector $e$ can be approximated in
a much lower-dimensional subspace, are the main ideas behind the multigrid 
methodology. Very efficient { geometric
multigrid} solvers have been developed for the case that the
lower-dimensional subspace corresponds to functions on a
well-structured grid (the {\it coarse} level), on which the smooth
errors can be approximated by easy-to-derive equations, based for
example on discretizing the same continuum operator that has given
rise to the fine-level equations $Ax=b$. The coarse-level
equations are solved  recursively using a similar combination of
relaxation sweeps and still-coarser-level approximations to the
resulting smooth errors.

To deal with more general problems, e.g., ones where the
fine-level system may not be defined on a well-structured grid
nor perhaps arise from any continuum problem, {algebraic
multigrid} (AMG) methods develop techniques for deriving the
set of coarse-level {\it variables} and the coarse-level {\it
equations}, based directly on the given matrix $A$. The basic
approach, developed in \cite{BMR83, oAMG, JWRuge_KStuben_1987a} and called today
classical AMG or RS-AMG,  involves the following two steps.
\bi
\item[(1)] Choosing the coarse-level variables as a subset $\coarsevar$ of
the set $\Omega$ of fine variables, such that each variable in $\Omega$ is
{\it strongly connected} to variables in $\coarsevar$.
 \item[(2)]
Approximating the fine-level residual problem $Ae=r$ by the
coarse-level equations $A_c e_c=r_c$ using the Galerkin
prescription $A_c = RAP$ and $r_c=Rr$,  yielding an approximation $ P e_c$ to $e$.
\ei 
The {\it interpolation matrix} $P$, and the {\it restriction matrix} $R$
are both defined directly in terms of the elements of the matrix
$A$, and this construction relies heavily on the notion 
of  {\it strong connections} developed for $M$-matrices.  

In the past two decades, many variations of the classical AMG algorithm have
been introduced, including modifications to the coarse-grid selection
algorithms
\cite{AJCleary_RDFalgout_VEHenson_JEJones_1998a,HDeSterck_UMYang_JJHeys_2005a,VEHenson_UMYang_2002a}, 
the definition of interpolation \cite{MBrezina_etal_2005a}, or both
\cite{Brannick_Trace_06,MBrezina_etal_2000a,TChartier_etal_2003a,PVanek_JMandel_MBrezina_1995a}.
A summary and comparison of many of these variations within the
classical AMG framework appears in \cite{GAMeurant_2001a}.  Such
variations arise because
 the applicability of the original algorithm and its variants
is limited by the $M$-matrix heuristics upon which they are based.
In particular, these variants of the classical AMG approach yield efficient solvers for problems where the (properly scaled) matrix $A$ has a dominant diagonal and, with small possible exceptions, all its off-diagonal elements have the same sign.  Even then, the produced interpolation can have limited accuracy, insufficient for full multigrid efficiency.  An example where 
performance of the classical AMG method can deteriorate, and the one we consider here, is given by scalar anisotropic diffusion problems. 

Recent advances in the development of  algebraic multigrid have focused on alternative notions of strength of connection that do not rely on specific properties of the system matrix, $A$, in their construction. 
  In \cite{EChow_2003a},  for example, 
samples of algebraically smooth error are used to determine the directions in
which this error varies slowly and, thus, interpolation can be very
effective.  For many problems, however, such a direction may not
exist, and the proposed technique does not apply.  An alternative approach 
was developed by Br\"{o}ker~\cite{OBroeker_2003a}, in which the relative sizes of the entries of the sparse approximate inverse preconditioning matrix \cite{MJGrote_THuckle_1997a} are used to
indicate strength of connection.  In some respects, this approach is a
special case of the approach proposed in~\cite{Bran_EBSOC},
which also measures strength using columns of the inverse; 
however, the measure used in
\cite{OBroeker_2003a} is based on the $L_2$ sizes of the approximate
inverse, not on the energy norm as  in~\cite{Bran_EBSOC}.
  Closely related to this work is the  approach of compatible relaxation
\cite{ABrandt_2000a,JBrannick_RFalgout,Brannick_Trace_06,PanayotRob_2003,OLivne_2004a} which 
uses a modified relaxation scheme to expose the character of the
slow-to-converge, i.e., algebraically smooth,  error.  Coarse-grid points are then selected where
this  error is the largest, thus avoiding explicit use of a measure of strength of connection
in choosing the coarse variables. 

In the present article, we focus on developing an alternative coarsening algorithm to select coarse variable and interpolatory sets for 
anisotropic systems where even these newly developed variants of AMG can often be difficult to apply.  
Our approach is based on a general notion of strength of connection derived from the minimum of the local least squares problem used in defining interpolation \cite{B00,BAMG2010}.  The basic approach finds interpolation to fit a set of test vectors (representatives of algebraically smooth error) in a least squares sense.  
We use the LS functional resulting from a simplified direct form of interpolation to construct an algebraic distance measure of strength of connection which is used to derive a strength
graph.   The strength graph is then passed to a coloring 
algorithm~\cite[Chap. 8]{WLBriggs_VEHenson_SFMcCormick_2000a} 
to coarsen the unknowns at each stage of a compatible relaxation 
coarsening algorithm.  
The quality of interpolatory sets are measured a posteriori, by first constructing LS interpolation to these sets and then monitoring the corresponding values of the LS functional for the computed minimizers \cite{B00,BAMG2010}.   
A main new feature  of our coarsening algorithm is that it allows for aggressive
and problem-oriented  coarsening as well as long-range interpolation in a straightforward way.  
This is accomplished by deriving long-range measures  of strength between neighboring pairs (or sets) 
of fine and coarse points.  

Several algorithms for aggressive coarsening and long-range interpolation have been developed in the past~\cite{KStuben_2001a, HDeSterck_UMYang_JJHeys_2005a}, yet these techniques rely on the same $M$-matrix based strength of connection measure used in the original classical AMG coarsening algorithm~\cite{BMR83, oAMG, JWRuge_KStuben_1987a}.    Aggressive coarsening is then achieved by applying the classical AMG coloring algorithm multiple times before constructing interpolation.   While such approaches have been shown to be effective in reducing the large complexities often resulting from application of the classical AMG scheme to finite-difference and finite-element discretizations of elliptic problems, they are again limited to $M$-matrix type problems.  In contrast, our proposed long-range measure aims to address the issue of strength of connection in a more general context -- the goal being to determine explicitly those degrees of freedom from which high quality least squares interpolation for some given set of test vectors can be constructed.  

The remainder of the paper is organized as follows.  First, in  Section \ref{sec:background}, we give an  overview of the classical algebraic multigrid coarsening approach.    Section \ref{sec:bamg} contains an introduction to the bootstrap algebraic multigrid components, with an emphasis on the least squares interpolation approach and compatible relaxation coarsening algorithm.   Then, in Section \ref{sec:distances}, a general definition of strength of connection and the notion of  algebraic distances,  as well as its connection to compatible relaxation, are discussed.   
The diffusion equation with anisotropic coefficients and its descritization are introduced in Section \ref{sec:numerics} as are  the results of numerical experiments of our approach applied to this system.
 
\section{Classical AMG Coarsening}
\label{sec:background}
The coarse-grid selection performed in the classical AMG algorithm
\cite{oAMG,JWRuge_KStuben_1987a} is made  based on  properties of   slow-to-converge errors discovered by applying  pointwise relaxation to  discrete systems with  diagonally dominant $M$-matrices.  
These  errors are  used to identify  the important (strong)
connections within the linear system.  The coarse-grid points are then
selected using maximal independent subset heuristics to ensure a
significant reduction in grid size, still maintaining accurate
approximation properties. The outline of the classical AMG coarsening procedure is given next.

Consider the Gauss-Seidel iteration, with error propagation operator
$I-L^{-1}A$, where  $L$ is the lower-triangular part of $A$.  Fixing this to
be the fine-scale relaxation used in a multigrid method for the
symmetric and diagonally dominant M-matrix, $A$, the goal of
coarse-grid correction is to effectively reduced  the error components
that are not significantly reduced by Gauss-Seidel relaxation.  It is
shown in \cite{oAMG} that such an error,
$e$, must satisfy
\[
\sum_{i,j} (-a_{ij})(e_i-e_j)^2 \ll \sum_i a_{ii}e_i^2.
\]
Thus, if $a_{ij}$ is large, relative to $\displaystyle \max_{k\ne i} |a_{ik}|$ or
$\displaystyle \max_{k\ne j} |a_{jk}|$, then it must be true that $e_i \approx
e_j$.

This observation leads to the definitions of strong dependence which
has significant influence on the selection of AMG coarse-grid points.  For a
given degree of freedom, $i$, the set of points that $i$ {\it strongly
depends upon} is defined as
\begin{equation}
S_i = \left\{ j : -a_{ij} \ge \theta \max_{k \ne i} 
\left\{-a_{ik}\right\} \right\},
\label{eq:classic_strength}
\end{equation}
for some suitable choice of $\theta$, $0 < \theta \le 1$.  Similarly, the set of points
that are {\it  strongly influenced by} $i$ is denoted $S_i^T$:  $j \in S_i^T$
implies that $i \in S_j$.

Once strong connections are determined, a coarse grid is chosen so
that all strongly connected neighbors of any fine-grid point (that is
not itself a coarse-grid point) are available for direct or a
path-length two indirect interpolation.  That is, if two fine-grid
points are strongly connected to one another, they must have a common
coarse-grid neighbor, so that this strong connection may be accounted
for indirectly in interpolation.  This condition is usually expressed
as the first AMG coarsening heuristic.\\

\begin{itemize}
\item[{\bf H1:}] For each fine-grid point $i$, each point $j \in S_i$ must
either be a coarse-grid neighbor of $i$ or strongly depend on at least one
coarse-grid neighbor of $i$.\\
\end{itemize}
If used alone, this rule tends to create large coarse grids, as it is
most easily satisfied by adding points to the coarse grid.  In
practice, a second heuristic is used to limit the size of the coarse
grid.\\

\begin{itemize}
\item[{\bf H2:}] The set of coarse-grid points, $C$, should be a maximal
subset of the original set of fine-grid points such that no
coarse-grid point strongly depends on any other coarse-grid point.\\
\end{itemize}

\noindent
Because these two goals may be contradictory, typical AMG implementations
rely on enforcing \Hone{}  and use \Htwo{} as a guide.  This is done by selecting 
coarse-grid points  using  a two-pass algorithm
that picks a set satisfying \Htwo{}, then checks for any points where
\Hone{} is violated, adding new coarse-grid points to compensate if
this occurs.  The first stage is often implemented as a coloring
algorithm \cite[Chap. 8]{WLBriggs_VEHenson_SFMcCormick_2000a}, where
coarse points are selected based on their number of strongly connected
neighbors.  Initially, all points are weighted with the number of
points that strongly depend upon them (that is, the size of $S_i^T$).
The point with the largest weight is then selected to be a coarse-grid
point.  Since each $j \in S_i^T$ is now strongly connected to a
coarse-grid point, all such $j$ are made fine-grid points so that
\Htwo{} is not violated.  All strongly connected neighbors of these
points (that is, $k\in S_j^T$ for any $j \in S_i^T$) are then made
more attractive as coarse-grid points, since they reflect unresolved
strong connections of fine-grid points.  Thus, the weights of all such
$k$ are incremented for each $j \in S_k$ that was made a fine-grid
point.  The algorithm then repeats, selecting the new largest-weighted
point as a coarse-grid point.

In this way, an initial coarse grid is chosen to give a maximal
independent set over all strong connections.  Then, 
if necessary, the minimal
number of points needed to ensure that property \Hone{} is fulfilled are added to the coarse grid. 
Once the selection of  coarse grid has been completed, an interpolation
operator is defined so that it is accurate for  errors that are slowly varying 
along strong connections.  We do not  provide specific details of AMG
interpolation here and, instead, again refer the reader to
\cite[Chap. 8]{WLBriggs_VEHenson_SFMcCormick_2000a}.  
We note, however, that the classical approach for constructing $P$ also 
relies on the above $M$-matrix heuristic that the direction in which the smooth error is locally constant
can be identified by the entries in the system matrix, thus, limiting its range of applicability.

\section{Bootstrap AMG}
\label{sec:bamg}
The bootstrap AMG (BAMG) algorithm, introduced in \cite[\S
17.2]{B00}, was developed to extend AMG to general (non $M$-matrix)
problems.
The BAMG approach combines the following two general devices to
inexpensively construct high-quality interpolation. \bi
\item[(A)] {The interpolation is derived to provide the best  least squares
fit  to a set of $\tau$-smooth 
test vectors (TVs) } obtained by a process described
below. \ei Denote by $C_i$ the set of coarse variables used in
interpolating to $i\in \finevar$. It follows from the
satisfaction of the compatible relaxation  criterion (see the next section) that, with a proper choice of
$C_i$, there exists an interpolation that will commit  only $\mathcal{O}(\tau)$
errors in reproducing $x_i$ for {\it any} vector $x$ which
is $\tau$-smooth. The size $|C_i|$ of this set should in principle
increase as $\tau$ decreases (and smaller $\tau$ means overall better multigrid convergence), but in practice a pre-chosen and sufficiently small interpolation {\it caliber} $\caliber := \max_{i \in \finevar}|\coarsevar_i|$,  
often yields small enough
errors. 
The set $C_i$ can often be adequately chosen by natural
considerations, such as the set of geometrical neighbors with
$i$ in their convex hull. If the chosen set is inadequate, the
least squares procedure will show {\it bad fitness} (interpolation
errors large compared with $\tau$), and the set can then be
improved. The least squares procedure can also detect variables in
$C_i$ that can be discarded without significant accuracy loss.
Thus, this approach allows creating interpolation with whatever
needed accuracy which is {\it as sparse as possible}. 
\bi \item[(B)] {Generally, the test vectors 
are constructed in a bootstrap manner}, in which several tentative
AMG levels are generated by interpolation fitted to only moderately smooth TVs; this tentative (multilevel)  structure  is then  used to produce {\it improved} (much
smoother) TVs, yielding a more accurate interpolation operator.  The process continues  if needed until fully efficient AMG
levels have been generated. \ei

The first test vectors are each produced by relaxing the homogeneous system
$Ax=0$ with a different starting vector. This quickly leads to a
$\tau$-smooth test set with $\tau \ll 1$. (A mixture of random
vectors and/or  diverse geometrically smooth vectors can generally be
used as initial  approximations. In the case of discretized isotropic
PDEs, if geometrically smooth vectors that satisfy the homogeneous
boundary conditions are used, relaxation may not be needed at all.
In many other cases relaxation can be confined to the neighborhood
of boundaries and discontinuities.) 
For the anisotropic problem considered, we use a constant vector and a set of, initially  positive random,    relaxed test vectors 
to define an algebraic distance notion of strength of connection and
the least squares interpolation operator.

\subsection{Least squares interpolation}\label{sec:lsp}

The basic idea of the least squares interpolation approach is to approximate
a set of test vectors, $\mathcal{V} = \{v^{(1)}, \ldots, v^{(k)}\}
\subset \mathbb{R}^n$, minimizing the interpolation error for these
vectors in a least squares sense.  In the context considered here,
namely, applying the least squares process to construct a classical AMG form of interpolation,  \textcolor{black}{each row of $P$, denoted by $p_i$,}
is defined as the minimizer of a local least squares functional. Given
a splitting of variables $\finevar = \Omega \setminus \coarsevar$ 
for each $i \in \finevar$ find $p_{i} \in \mathbb{R}^{n_{c}}, n_c = 
|\coarsevar|$, such that
\begin{equation}\label{eq:LSfuncrowi}
\mathcal{LS}(p_i) = 
\sum_{\kappa=1}^k\omega_\kappa\left(v_{i}^{(\kappa)} - \sum_{j\in \coarsevar_i} \left(p_{i}\right)_{j} v_j^{(\kappa)}\right)^{2} \mapsto \min,
\end{equation} 
where $\coarsevar_{i} \subset C$.
The weights $\omega_{\kappa}>0$ can be chosen to reflect the global algebraic smoothness
of the test vectors.  We give our specific choice in the numerical experiments section. 

\begin{remark}
We note that if the restricted vectors $v^{(1)}_{\coarsevar_i},\ldots,v^{(k)}_{\coarsevar_i}$ form a
basis for the local linear space $\mathbb{R}^{n_i}$, $n_i = |C_i|$,
then the solutions to the local least squares minimization problems are unique.
This in turn suggests setting a lower bound on
the number of vectors, $k$, used in the least squares fit
\begin{equation*}
  k \geq  \caliber.
\end{equation*}
Further, in practice we have observed that the accuracy of the least
squares interpolation operator and, hence, the performance of the resulting solver generally improves with increasing $k$~\cite{BAMG2010}, up to some value proportional to caliber $\caliber$. 
\end{remark}

\begin{remark}
In~\cite{KKahl_2009}, it was observed that the implicit application
of a local Jacobi relaxation to the TVs used in the least squares definition of interpolation is equivalent to an operator induced
form of least squares interpolation constructed using an Element-free AMG type approach (see~\cite{Vassilevski_2005}).  
\textcolor{black}{This equivalence of the least squares approach with an additional local relaxation step was also formulated and discussed in a slightly different scope in~\cite{iBAMG}, 
where it was defined as a residual-based least squares fit
\begin{equation}\label{eq:modls}
  \mathcal{LS}(p_i) = 
\sum_{\kappa = 1}^k\omega_{\kappa} \left(v_{i}^{(\kappa)} - \frac{1}{a_{ii}}r^{(\kappa)}_{i} - \sum_{j\in \coarsevar_{i}} (p_{i})_j v_{j}^{(\kappa)}\right)^{2} \mapsto \min, 
\end{equation}
which in turn was shown to be consistent with a classical AMG operator
induced form of least squares interpolation.} 
\end{remark}


\subsection{Compatible relaxation}\label{sec:CR}
A general criterion for choosing an adequate set of coarse
variables is the fast convergence of  {\sl compatible relaxation}
(CR), as introduced in \cite{ABrandt_2000a} (being a special case of
choosing coarse variables for much more general types of systems
\cite{SU}, introduced first for problems in statistical mechanics
\cite{RMG}). In fact, an improved version called
habituated compatible relaxation, introduced in \cite{OLivne_2004a}, yields an
accurate prediction for the convergence rate that can be achieved
by an AMG solver that employs the given relaxation scheme and the
proposed set of coarse variables. Analysis of general compatible relaxation  approaches
and works on the development of compatible relaxation  coarsening algorithms are found 
in~\cite{PanayotRob_2003,Brannick_Trace_06,JBrannick_RFalgout}.  We note that the ability of habituated compatible relaxation 
to quantitatively predict performance of the solver
is very useful in designing and debugging the actual solver.


Given  the coarse-and-fine level splitting, $\coarsevar$ and $\finevar$, 
a simple form of 
compatible relaxation is given by $\finevar$-relaxation for the homogeneous system --- relaxation applied only to the set 
of $\finevar$ variables. 
Given the partitioning of $\Omega$ into $\finevar$ and $\coarsevar$,
we have
$$
u = \begin{pmatrix} u_f \\ u_c \end{pmatrix}, \enspace 
B = 
\begin{pmatrix}
    B_{ff}  & B_{fc}   \\
    B_{cf} & B_{cc} 
\end{pmatrix},
\enspace \mbox{and} \enspace
M = 
\begin{pmatrix}
    M_{ff}  & M_{fc}   \\
    M_{cf} & M_{cc} 
\end{pmatrix},
$$
assuming the equations are permuted such that the unknowns in $\finevar$ 
come before those in $\coarsevar$. 
The $\finevar$-relaxation form of compatible relaxation is then defined by 
\begin{equation} \label{frelax:eq}
   u_f^{\nu+1} = (I-M_{ff}^{-1}B_{ff})u_f^{\nu} = E_{f}u_{f}^{\nu}. 
\end{equation}
If $M$ is symmetric, the asymptotic convergence rate of compatible relaxation,   
\[
 \rho_{f}
=\rho(E_{f}), 
\]
where $\rho$ denotes the spectral radius, provides a measure of the quality of the
coarse space, that is, a measure of the ability of the set of coarse
variables to represent error not eliminated by the given fine-level relaxation.
This measure can be approximated using $\finevar$-relaxation for the homogeneous 
system with a random initial guess $u_f^{0}$. Since $\lim_{\nu \to \infty} 
\|E_f^\nu\|^{1/\nu} = \rho(E_f)$ for any norm $\| \cdot \|$, the measure
\begin{equation} \label{meas1:eq}
\varrho_f  = \left(\|u_f^\nu\| / \|u_f^0\| \right)^{1/\nu}
\end{equation}
estimates $\rho_f$. 

The set of coarse
variables is then constructed using a multistage coarsening algorithm,
where a single stage consists of: (1) running several iterations of compatible relaxation 
(based on the current set $F$) and (2) if it  is slow to converge,
adding an independent set of the fine-level variables not well treated by CR to $C$.  
Steps (1) and (2) are applied repeatedly, giving rise to a sequence of
coarse variable sets:
\[\emptyset = C_0 \subseteq C_1 \subseteq ... \subseteq C_m, \]
until convergence of compatible relaxation   below a 
prescribed tolerance is reached, yielding an accepted coarse set $C := C_m$. 

 We note that all current CR algorithms have intentionally avoided the explicit use of 
a strength of connection measure in constructing the coarse variable sets, instead relying on the error
produced by CR to form candidate sets of potential $\coarsevar$-points.  Our aim here is to develop
a more general notion of strength of connections based on algebraic distances 
and to explore its use in a compatible relaxation coarsening scheme
and defining AMG interpolation.
  
\section{Selecting coarse variables and interpolation via algebraic distances}
\label{sec:distances}
While the classical definition of strong dependence is appropriate for the case
of diagonally dominant M-matrices for which it was intended, it is
frequently seen to break down when applied to problems involving other classes of 
matrices.
The near null space of a diagonally dominant M-matrix is typically a locally slowly varying (or locally
constant) vector, and, in such cases the AMG strength heuristic relies on this being
reflected in the coefficients of the system matrix, $A$.  If either the near
null space cannot be accurately characterized as locally constant or
this is not reflected in the matrix coefficients, then AMG performance
typically suffers.

To derive a more general measure of strength of connection, we 
consider qualifying strength among variables based on a variables ability to
interpolate $\tau$-smooth error to its neighbors.  We thus characterize a given variable (set of variables) as strongly connected
to a neighboring grid point as one(s) for which the value of the LS functional is small when building interpolation
from this point (set of points), for some set of test vectors.  In doing so, we are able to identify
suitable coarse-grid points as those from which it is possible to build a high-quality LS interpolation
to its neighbors.  Next, we describe the idea of  measuring strength between neighbors via algebraic distances in detail and then discuss how we use this measure in a CR coarsening algorithm and in computing interpolation.

\subsection{Strength of connection by algebraic distances}\label{sec:alg_dist}

The basic idea of the {\em algebraic distances} approach to measuring strength of connection is to construct a row of least squares interpolation, $p_i$, for each fine variable $i\in \Omega$ and various choices of its interpolatory set, $C_i$:
\begin{equation}\label{eq:LSfunc}
  \mathcal{LS}(p_i) = 
\sum_{\kappa = 1}^k\omega_{\kappa} \left(v_{i}^{(\kappa)} - \frac{1}{a_{ii}}r^{(\kappa)}_{i} - \sum_{j\in \coarsevar_{i}} (p_{i})_j v_{j}^{(\kappa)}\right)^{2} \mapsto \min , 
\end{equation}
and then using the associated values of $\mathcal{LS}$ to define weights quantifying the connectivity among variables.   

With a given matrix $A$, associate a connected
graph $\mathcal{G}=(\mathcal{V},\mathcal{E})$,
where $\mathcal{V}$ and $\mathcal{E}$ are the sets of vertices and
edges, with $n=|\mathcal{V}|$ (cardinality of $\mathcal{V}$).  Here, an edge $(i,j) \in E \iff (A)_{ij} \neq0$.   
Let $\mathcal{G}_d(\mathcal{V}_d,\mathcal{E}_d)$ denote the graph of the matrix $A^d$ and
define $\mathcal{G}_{d, i}(\mathcal{V}_i,\mathcal{E}_i)$ as  the subgraph associated with the $i$th vertex and its algebraic neighbors:
\begin{equation}\label{graph}
\mathcal{V}_{i} := \{\:  j \:  | \: (A^d)_{ij} \neq 0\} \quad \mbox{and} \quad \mathcal{E}_i := \{ (i,j) |  \: (A^d)_{ij} \neq 0\}.
\end{equation}  

In its simplest form, the notion of algebraic distances is
straightforward: For a given search depth, $d$, set of test vectors, $\{v^{(1)},  ..., v^{(k)}\}$, and fine variable $i\in \Omega$ compute
\begin{equation}
 r_{ij}  = \frac{1}{
\sum_{\kappa = 1}^k\omega_{\kappa} \left(v_{i}^{(\kappa)} - \frac{1}{a_{ii}}r^{(\kappa)}_{i} -  p_{ij} v_{j}^{(\kappa)}\right)^{2}} \quad \mbox{for all} \quad j \in \mathcal{V}_i , 
\label{eq:measure}
\end{equation}
where $p_{ij}$ is the minimizer of the least squares problem to a
single variable $j$.
Here, although the measure is able to determine the coupling between any given two points, we limit its use to local neighborhoods, $j \in \mathcal{V}_i$.  This simplification coupled with the idea of deriving strength according to an algebraic distance measure based on simple one-point (one-sided) interpolation, allows us to control the complexity of the algorithm.
More generally, the algebraic distance measure can be defined from the LS functional obtained by building interpolation from sets of points, thereby providing a local a posteriori measure of the quality of the interpolatory set $\coarsevar_i$.  

\begin{remark}
Using the given weighted strength graph in a coloring algorithm to aggregate unknowns or in a greedy algorithm to define a matching in the graph results in a so-called plain aggregation scheme, which in turn can be used as a solver within an AMLI cycle~\cite{brannick_amli_2011}.   In~\cite{brannick_markov_2011}, such a scheme was developed for finding steady state vectors of Markov chains.  A similar technique, which chooses aggregates using a greedy strategy based on a local stability measure is developed and analyzed in~\cite{brannick_local_stab_2011}.
\end{remark}

\subsection{Compatible relaxation coarsening and algebraic distances}\label{sec:cr}
In choosing $\coarsevar$, we integrate the simplified variant of the algebraic distance notion of strength of connection based on one-sided interpolation~\eqref{eq:measure} into the CR-based coarsening 
algorithm developed in~\cite{JBrannick_RFalgout}.   The notion of algebraic distances is used to form a subgraph of the graph of  
the matrix $A^d$, $d=1,2,\dots$
Specifically, the algebraic distance between any two adjacent vertices 
in the graph, $\mathcal{G}_d$, 
of $A^d$ is computed using~\eqref{eq:measure}.  Then, for each vertex, $i \in \finevar$, we remove edges
adjacent to $i$ with small weights relative to the largest weight of all edges connected to $i$:
\begin{equation}\label{eq:strength_graph}
\mathcal{V}_{\mathcal{M}} = \finevar \; ; \quad 
\mathcal{E}_{\mathcal{M}} = \{ (i,j) \:  | \: i,j \in \finevar  \quad \mbox{and} \quad r_{ij} > \theta_{ad}  \max_k r_{ik} \}, 
\end{equation}
with $ \theta_{ad} \in (0,1)$.
This in turn
produces the graph of strongly connected vertices, $\mathcal{M}^d(\mathcal{V}_{\mathcal{M}}, \mathcal{E}_{\mathcal{M}}) $.  
Note that by definition, the strength graph is restricted to vertices $i\in F$ so that it can be used in subsequent CR coarsening stages.

The strength graph, $\mathcal{M}^d$, is then passed to a coloring
algorithm, as in the classical AMG approach outlined in Section~\ref{sec:background}, where
as mentioned earlier, coarse points are selected based on their number of strongly connected
neighbors.  Recall that, in the classical approach strong couplings are defined in terms of $S_i$ and $S_i^T$ as in Equation \eqref{eq:classic_strength}.  In our approach,  the sets $S_i$ and $S_i^T$ containing the information on strong coupling are instead derived from the algebraic distance based connectivity graph $\mathcal{M}^d$.

A rough description of our overall CR coarsening approach is described by Algorithm 1;  for additional specific details of the CR algorithm we refer the reader to~\cite{JBrannick_RFalgout}, in which this scheme was developed for diffusion problems similar to those we consider here.

\begin{algorithm}
\caption{compatible\_relaxation \hfill \{Computes $\coarsevar$ using Compatible Relaxation\}\label{alg:CR}}
\begin{algorithmic}
\STATE \textit{Input:} $\coarsevar_{0}$ \qquad \COMMENT{$\coarsevar_{0} = \emptyset$ is allowed}.
\STATE \textit{Output:} $\coarsevar$
\STATE Initialize $\coarsevar = \coarsevar_{0}$
\STATE Initialize $\mathcal{Z} = \Omega \setminus \coarsevar$
\STATE Perform $\nu$ CR iterations with constant $u^0$
\WHILE{$\rho_{f} > \delta$}
\STATE $\mathcal{Z} = \{i \in \Omega \setminus \coarsevar: 
\sigma_i  > \mathbf{tol} \}$
\STATE $\coarsevar = \coarsevar \cup \{\text{\: independent set of} \: \mathcal{Z} \: \}$ based on $\mathcal{M}^d$
\STATE Perform $\nu$ CR iterations with constant $u^0$
\ENDWHILE
\end{algorithmic}
\end{algorithm}

\begin{remark}
Various choices of the candidate set measure, $\sigma_i$, used in determining potential $\coarsevar$-points have been studied in the literature~\cite{ABrandt_2000a,OLivne_2004a,JBrannick_RFalgout}; our's follows the  definition in~\cite{JBrannick_RFalgout}, namely 
$$ \sigma_i := \frac{\displaystyle |{u}_{i}^{\nu}|}{\displaystyle \|{u}^{\nu}\|_{\infty}}.$$  
We mention that, in practice, this measure gives the best overall results for smaller values of $\nu$, say $\nu=5$.  
Additional discussions and extensive testing of this approach are found in~\cite{JBrannick_RFalgout}.
\end{remark}

\begin{remark}
The matrix $A^d$, $d>1$, is never formed in practice --  it is straightforward to reconstruct the connectivity graph of $A^d$ using local operations involving only the graph of $A$.
\end{remark}

\subsection{Defining interpolation via algebraic distances}\label{sec:greedy}
Given a set of coarse variables $\coarsevar$ and a set of $\tau$-smoooth test
vectors, $\{v^{(1)},  ..., v^{(k)}\}$, interpolation of a given
caliber, $c$, is constructed
via the least squares functional defined in~\eqref{eq:LSfunc}. 
More specifically, we consider all possible 
sets of $\coarsevar$-points, $W$, up to a given cardinality as prescribed by
parameter, $c$,  in the $d_{LS}$-ring coarse point neighborhood of a
given $\finevar$-point, $i$,  defined as 
\begin{equation}\label{eq:nbhd}
\mathcal{N}_{d_{LS}, i} :=  \coarsevar \cap \mathcal{V}_{d_{LS}, i}, 
\end{equation}
where the connectivity in $\mathcal{G}_{d_{LS}, i}$  is defined as in~\eqref{graph}.  
Thus, the sets of possible interpolatory points can be written as 
$$\mathcal{W} := \{ W | \;  \: W \subseteq \mathcal{N}_{d_{LS}, i} \; \mbox{and} \;  |W| \leq c  \}.
$$ By sampling this set we find the minimizer
 of the least squares functional: for each $i \in \finevar$
$$
C_i = \arg\min_{W \in \mathcal{W}} \mathcal{LS}_W(p_i) , $$
where, 
$$
 \mathcal{LS}_W(p_i) = \sum_{\kappa = 1}^k\omega_{\kappa} \bigg(v_{i}^{(\kappa)} - \frac{1}{a_{ii}}r^{(\kappa)}_{i} - \sum_{j\in W} (p_{i})_j v_{j}^{(\kappa)}\bigg)^{2}, 
$$
 and $p_i$ denotes the minimizer of the least squares problem \eqref{eq:modls} for the given set $W$.  Thus, we must compute $p_i$ and evaluate $\mathcal{LS}_W(p_i)$ for all possible choices of $W \in \mathcal{W}$.  The row of interpolation, $p_i$, is then chosen as the one that minimizes $\mathcal{LS}_W(p_i)$.

An additional detail of our approach is the penelization of large
interpolatory sets. It is easily shown that for two sets $W^{\prime}
\subset W^{\prime \prime}$, their corresponding minimal least squares
functional values fulfill $\mathcal{LS}_{W^{\prime}} \geq
\mathcal{LS}_{W^{\prime \prime}}$. In order to keep the interpolation
operator as sparse as possible, we require that the least squares
functional is reduced by a certain factor when increasing the
cardinality of $W$.
That is, a new set of points ${W^{\prime \prime}}$ is preferred
over a set of points ${W^{\prime}}$ with $|W^{\prime \prime}|
> |W^{\prime}|$ if
\begin{equation*}
  \mathcal{LS}_{W^{\prime \prime}} < \left(\mathcal{LS}_{W^{\prime}}\right)^{\gamma \left( |W^{\prime\prime}|-|W^{\prime}|\right)}.
\end{equation*} Based on numerical experience we usually choose $\gamma
= 1.5$  which  tends to produce accurate and sparse interpolation
operators for a large class of problems.

\begin{remark}
The above sampling of all possible combinations of interpolatory sets with cardinality up to some given caliber is one of many possible strategies for selecting $\coarsevar_i$.  In our experience, such an exhaustive search is rarely necessary and in most cases scanning a small number of  possibilities using some systematic strategy that incorporates the algebraic distance strength measure based on one-sided interpolation is sufficient.  
\end{remark}

\section{Numerical Results}
\label{sec:numerics}
In this section, to demonstrate the effectiveness of the algebraic distances measure of strength of connection when combined with CR coarsening and a greedy approach for constructing LS interpolation, we present tests of the approach applied to a variety of 2D anisotropic diffusion problems discretized on a $(N+1)\times(N+1)$ uniform grid.

\subsection{Model problem and its discretization}
We consider a finite-difference discretization of the two-dimensional diffusion operator 
\begin{equation}\label{eq:diff} \mathcal{L}\, u =  \nabla \cdot \mathcal{K} \nabla u 
\end{equation}
with anisotropic diffusion coefficient 
\begin{equation}  \mathcal{K}  =       \biggl( \begin{matrix}   \cos \alpha & -\sin \alpha \\
      \sin \alpha & \cos \alpha\ \\
   \end{matrix}\biggr)     \biggl( \begin{matrix} 
         1 & 0 \\
         0 & \epsilon \\
       \end{matrix} \biggr)   \biggl( \begin{matrix}   \cos \alpha & \sin \alpha \\
      -\sin \alpha & \cos \alpha\ \\
   \end{matrix}\biggr), 
   \label{eq:a}
   \end{equation} where
$0 <\epsilon<1$ and $0 \leq \alpha  < 2 \pi$. Changing   variables
\begin{equation}   \biggl( \begin{matrix}       \xi \\
      \eta  \\
   \end{matrix} \biggr)
   = \biggl( \begin{matrix}   \cos \alpha & \sin \alpha \\
      -\sin \alpha & \cos \alpha\ \\
   \end{matrix}\biggr) \biggl( \begin{matrix}       x \\
      y  \\
   \end{matrix} \biggr),
   \label{eq:change_var}
   \end{equation}
yields  strong connections aligned with direction $\xi$:
  \begin{equation} \mathcal{L} \, u(\xi, \eta) = u_{\xi \xi} + \epsilon u_{\eta \eta}.
  \end{equation}
Its  equivalent     
formulation in $(x,y)$ coordinates  is given by 
\begin{eqnarray}\label{eq:scal} \mathcal{L}u(x, y) = a u_{xx} + b u_{xy} + c u_{yy},  \end{eqnarray}
where $a  = \cos^2 \alpha + \epsilon \sin^2 \alpha, \, b = (1-\epsilon) \, \sin 2\alpha $, and $c = \sin^2 \alpha + \epsilon \cos^2 \alpha$.

In formulating a finite-difference  discretization of (\ref{eq:scal}), we consider a standard five-point discretization for the Laplacian term and then define the discretization of the $u_{xy}$ term using intrinsic  strength of connections.  In the $\alpha = \pi/4$ case, for example, this amounts to using lower-left  and upper-right  neighbors and avoiding  lower-right  and upper-left  ones  in defining the discretization for each fine-grid point $i \in \Omega$. The stencil for $u_{xy}$ and this choice of $\alpha$ is then given by 
\begin{equation}
\label{fd_seven} \tilde{S}_{xy} = \dfrac{1}{2h^2}  
\Biggl(
\begin{matrix}
0 & - 1 & 1 \\
0 & 1 & - 1 \\
0 & 0 & 0 \\
\end{matrix}
\Biggr) + \dfrac{1}{2h^2}  \Biggl(
\begin{matrix}
0 & 0 & 0 \\
-1 & 1 & 0 \\
1 & -1 & 0 \\
\end{matrix}
\Biggr).               
\end{equation}        
The overall stencil thus includes seven nonzero values.   In our numerical experiments, we consider also the worst case scenario, that is, the case where the discretization for $u_{xy}$ is defined  along  the weakest connection, $\eta$.   For instance, the appropriate stencil for $\alpha = \pi/4$ is a poor choice for $\alpha = -\pi/4$ as directions $\xi$ and $\eta$, as in (\ref{eq:change_var}), interchange for these choices of $\alpha$. Further, for this value of $\alpha$ and taking $\epsilon = .1$, the resulting system matrix has stencil
$$
S_A =  \Biggl(
\begin{matrix}
 & -1 & .45 \\
-1 & 3.1 & -1 \\
.45 & -1 &  \\
\end{matrix}
\Biggr),
$$ and, thus, is not an $M$-matrix.  Here, some of the off-diagonal entries in $A$ are positive and, hence, the heuristics motivating the classical definition of strength of connection are not applicable.  We consider this extreme case, although it is unlikely to arise in practice, to demonstrate the robustness  of our approach in  choosing the right coarsening strategy for the targeted anisotropic problems.   
We mention that on a structured grid, our chosen seven-point discretization is equivalent to the finite element discretization of the same
elliptic boundary value problem.
%

\subsection{Formulating a two-level coarsening algorithm}
Our aim is to study the robustness of the algebraic distance notion of
 strength of connection for grid and non-grid aligned anisotropy.  We focus our tests on 
 the seven-point finite-difference discretization introduced earlier in this section for various values
 of anisotropy angle, $\alpha$, anisotropy coefficient, $\epsilon$.
 
In all tests, coarse grids are chosen using the compatible relaxation approach discussed in Section \ref{sec:cr} and then interpolation is computed using the greedy LS approach given in Section \ref{sec:greedy}.  
The stopping tolerance for CR steps is set as $ \delta = 0.7$ and the number of CR sweeps is set as $\nu = 5$, so that the algorithm terminates when $\rho_f$ in \eqref{meas1:eq} is below this value, and the tolerance for the candidate set measure is taken as ${ \bf tol} = 1-\rho_f$.  The larger choice of stopping tolerance allows more aggressive coarsening whereas our choice of $\bf{tol}$ ensures that points are added to the candidate set more sparingly at later stages of the CR coarsening scheme.
We fix the strength of connection parameter used in forming the strength graph at $\theta_{ad} = .5$ and vary the the graph distance, $d$, used to define the graph $\mathcal{G}_d$, from which $\mathcal{M}_d$ in \eqref{eq:strength_graph} is constructed.  We note that generally the overall quality of the grids the algorithm produces depends only mildly on the choice of $\theta_{ad}$.

The search depth, used in defining the greedy approach for choosing LS interpolation as in \eqref{eq:nbhd}, is fixed as $d_{LS} = d + 2$.  Here, taking a larger value of the search 
depth than the coarsening depth allows the approach to scan a larger number of possible interpolatory sets in forming long-range interpolation (whenever the problem requires it).  In this way, the LS scheme for constructing interpolation is able to better follow a wider range of values of the anisotropy angle, $\alpha$. 
Eight test vectors are used to construct the modified LS interpolation given in \eqref{eq:modls},  seven are computed by applying 40 iterations of Gauss Seidel to distinct random initial guesses and the eighth is the constant vector also relaxed with 40 iterations of Gauss Seidel.  
We mention that the number of test vectors and especially the amount
of relaxation applied to them can be sufficiently reduced by replacing
a one-grid procedure by a multilevel bootstrap cycling scheme~\cite{BAMG2010}.
However, the main focus of this paper is the study of the performance of the algebraic
distances measure of strength of connections for the targeted anisotropic problems; 
the development of a multilevel algorithm is subject of 
on-going research.

When presenting results of the two-grid solver constructed by our algorithm we use two pre- and post- Gauss Seidel relaxation sweeps on the fine grid and a direct solver on the coarse grid.  Here, we use two pre- and post- smoothing steps because our algorithm generally produces aggressive coarsening.  The estimates of the asymptotic convergence rates are computed as follows:
\begin{equation*}
\rho = \frac{\|e^{\eta}\|_A}{\|e^{\eta-1}\|_A},
\end{equation*}
where $e^\eta$ is the error after $\eta=100$ MG iterations, applied to the homogenous system 
starting with random initial approximations.  We also report the operator complexity ratio
 $$\gamma_o =  \dfrac{nnz(A) + nnz(A_c)}{nnz(A)},$$
 and the coarsening factor
 $\gamma_g =  \dfrac{|C|}{|\Omega|}$.

\subsection{Choosing coarse variables and interpolation}
\begin{figure}
  \subfigure[$\alpha = 0, \rho = .02, \rho_f = 0, \gamma_o = 1.89, \gamma_g = .516$ \label{fig:s10}]{\includegraphics[scale = 0.5]{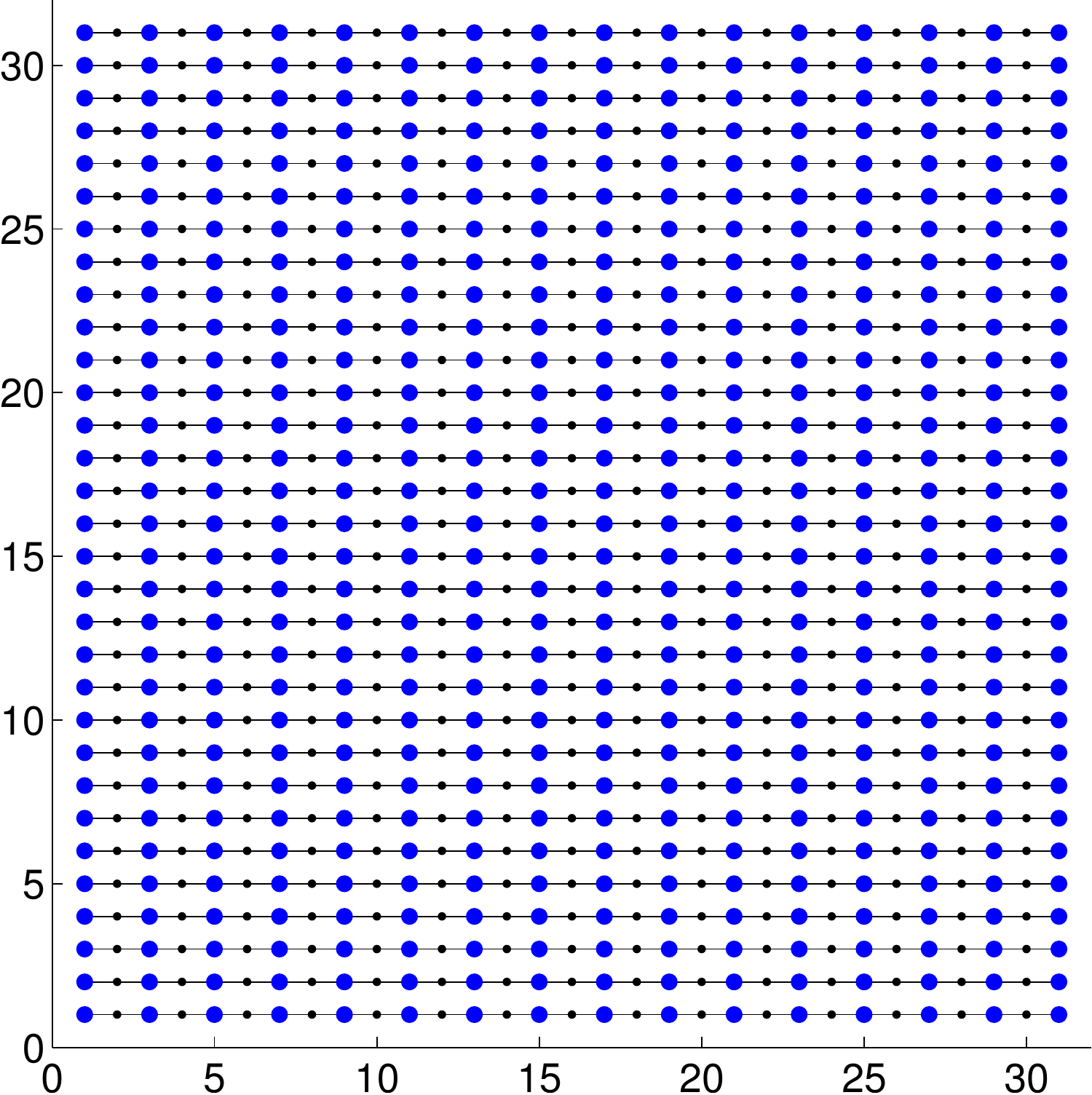}}
  \hfill
  \subfigure[$\alpha = \pi/4, \rho = .07, \rho_f = .50, \gamma_o =1.612, \gamma_g = .498$\label{fig:s1p4}]{\includegraphics[scale = 0.5]{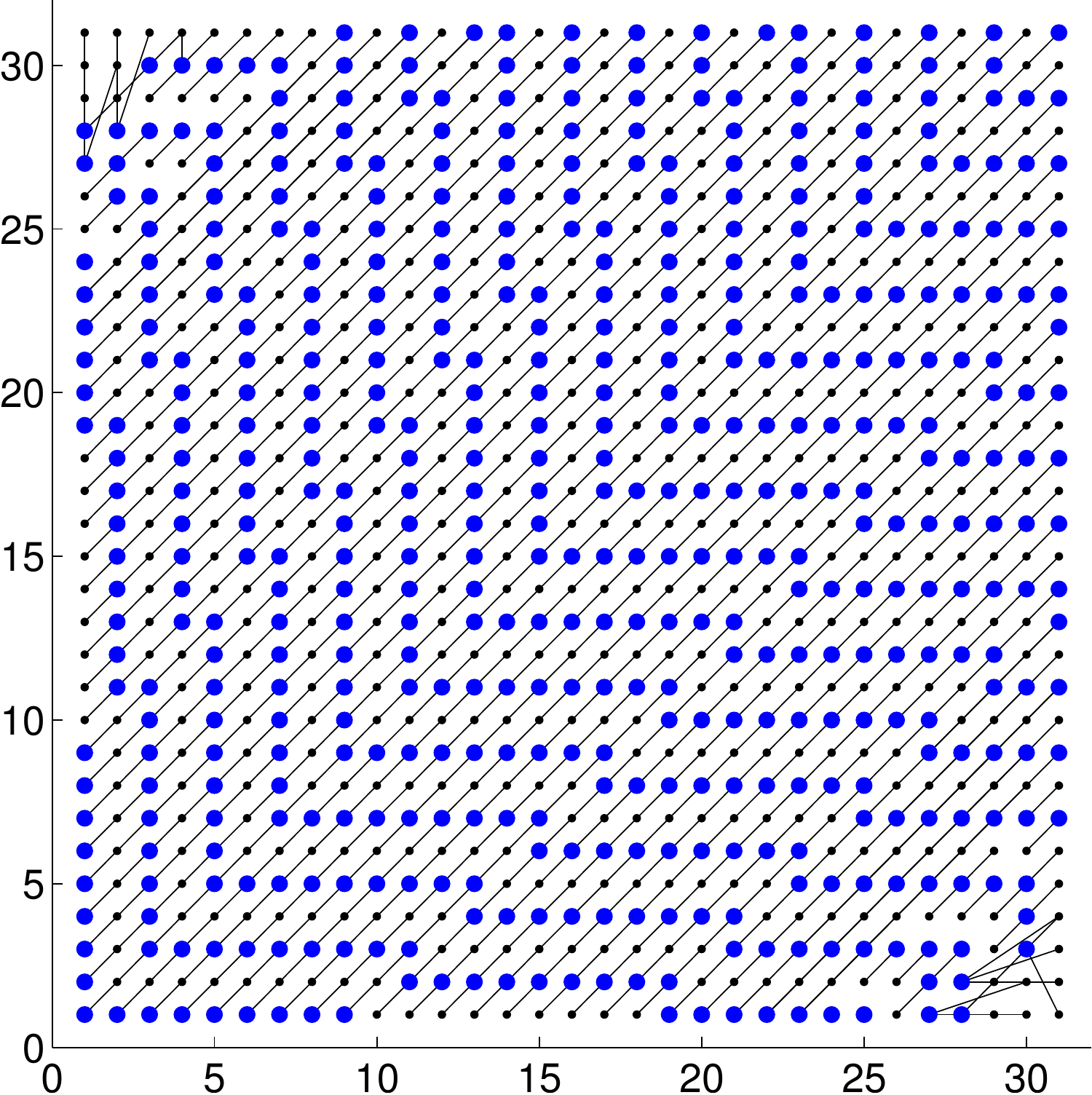}}
\hfill
  \subfigure[$\alpha = -\pi/4, \rho = .46, \rho_f = .78, \gamma_o = 1.556, \gamma_g = .354 $ \label{fig:s1mp4}]{\includegraphics[scale = 0.5]{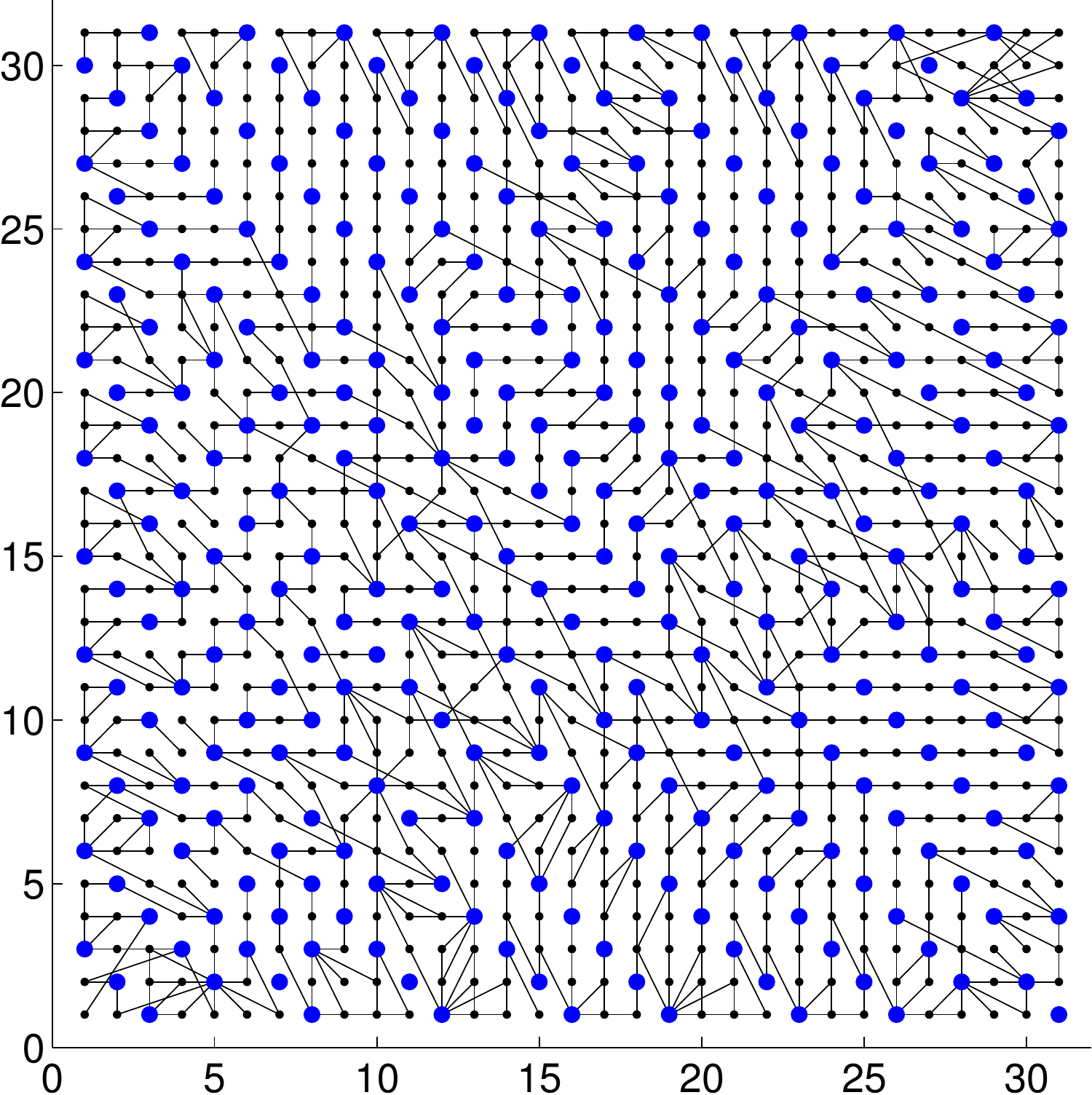}}
 \hfill
    \subfigure[$\alpha = \pi/8,   \rho = 0.45, \rho_f = 0.62, \gamma_o  = 1.910, 
\gamma_g = 0.426
 $\label{fig:s1p16}]{\includegraphics[scale = 0.5]{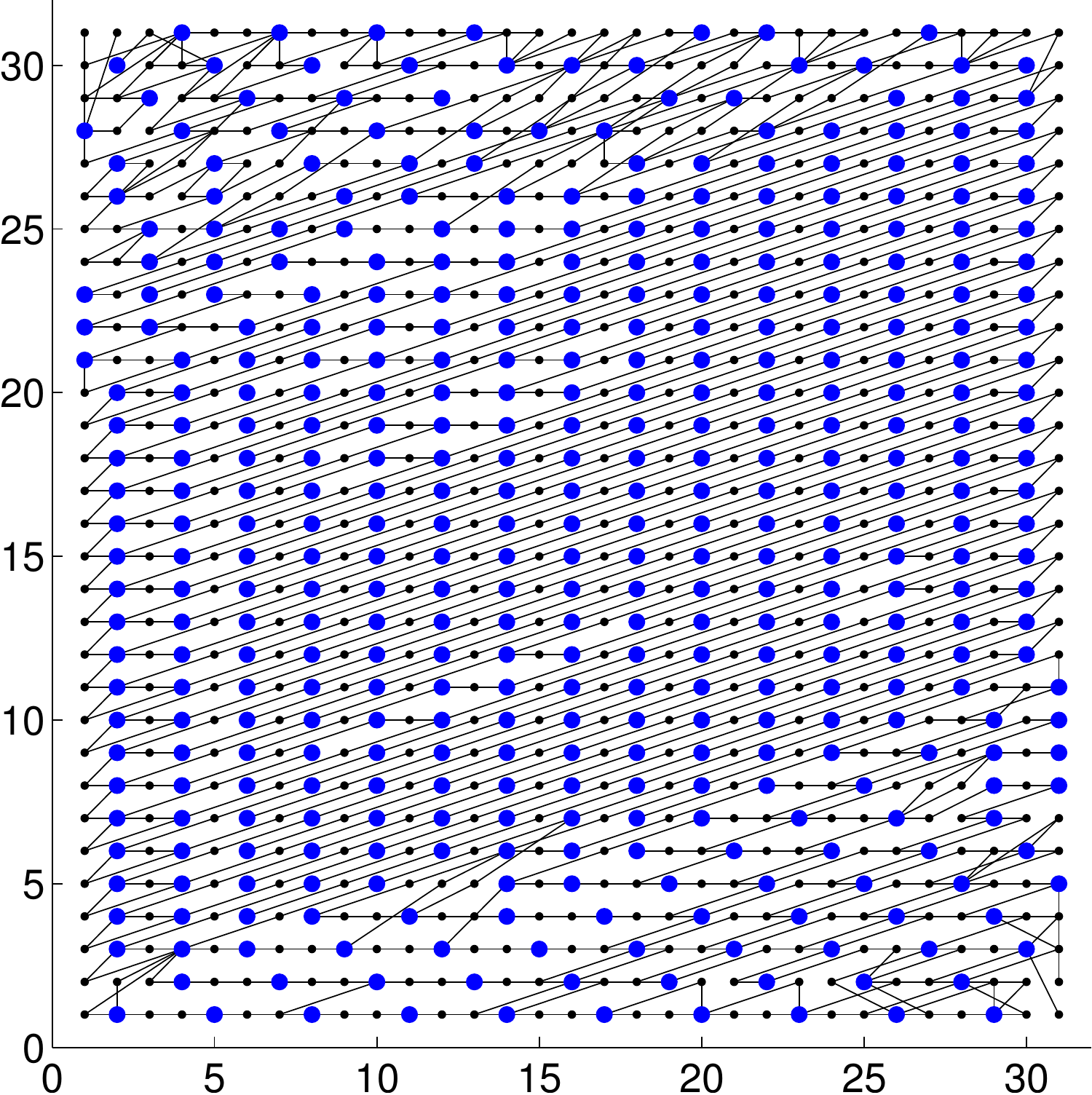}}
    \caption{Coarse grids and interpolation patterns for $h= 1/32$ for various choices of $\alpha$, using the graph of $A$, i.e., $d = 1$ and $d_{LS} = 3$,   to define the strength matrix. \label{fig:s1}}
\end{figure}

\begin{figure}
  \subfigure[$\alpha = 0, \rho = .05, \rho_f = .25, \gamma_o = 1.568, \gamma_g = .334 $\label{fig:s20}]{\includegraphics[scale = 0.5]{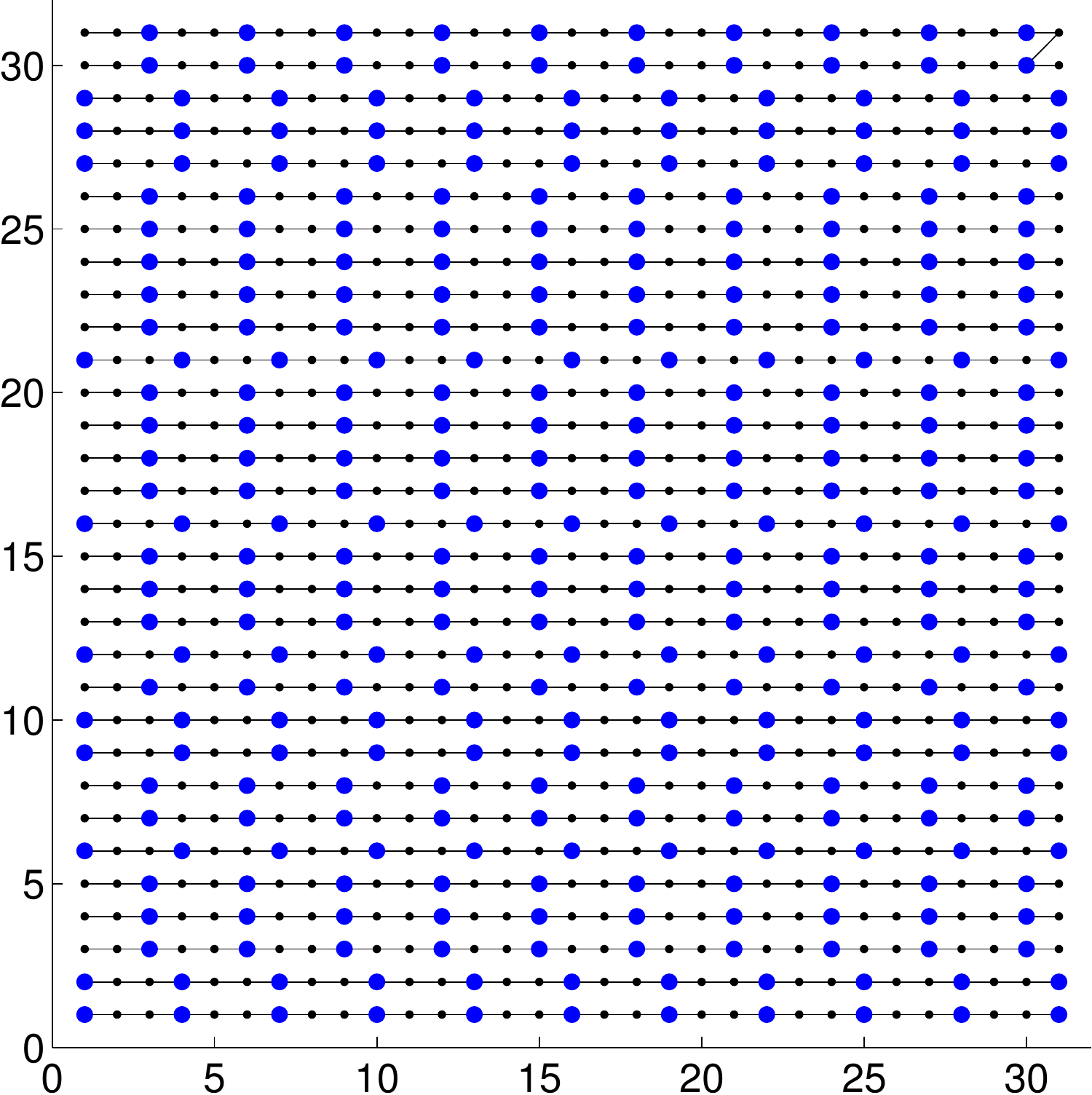}}  \hfill
  \subfigure[$\alpha = \pi/4, \rho = .06, \rho_f = .50, \gamma_o = 1.398, \gamma_g = 0.331 $\label{fig:s2p4}]{\includegraphics[scale = 0.5]{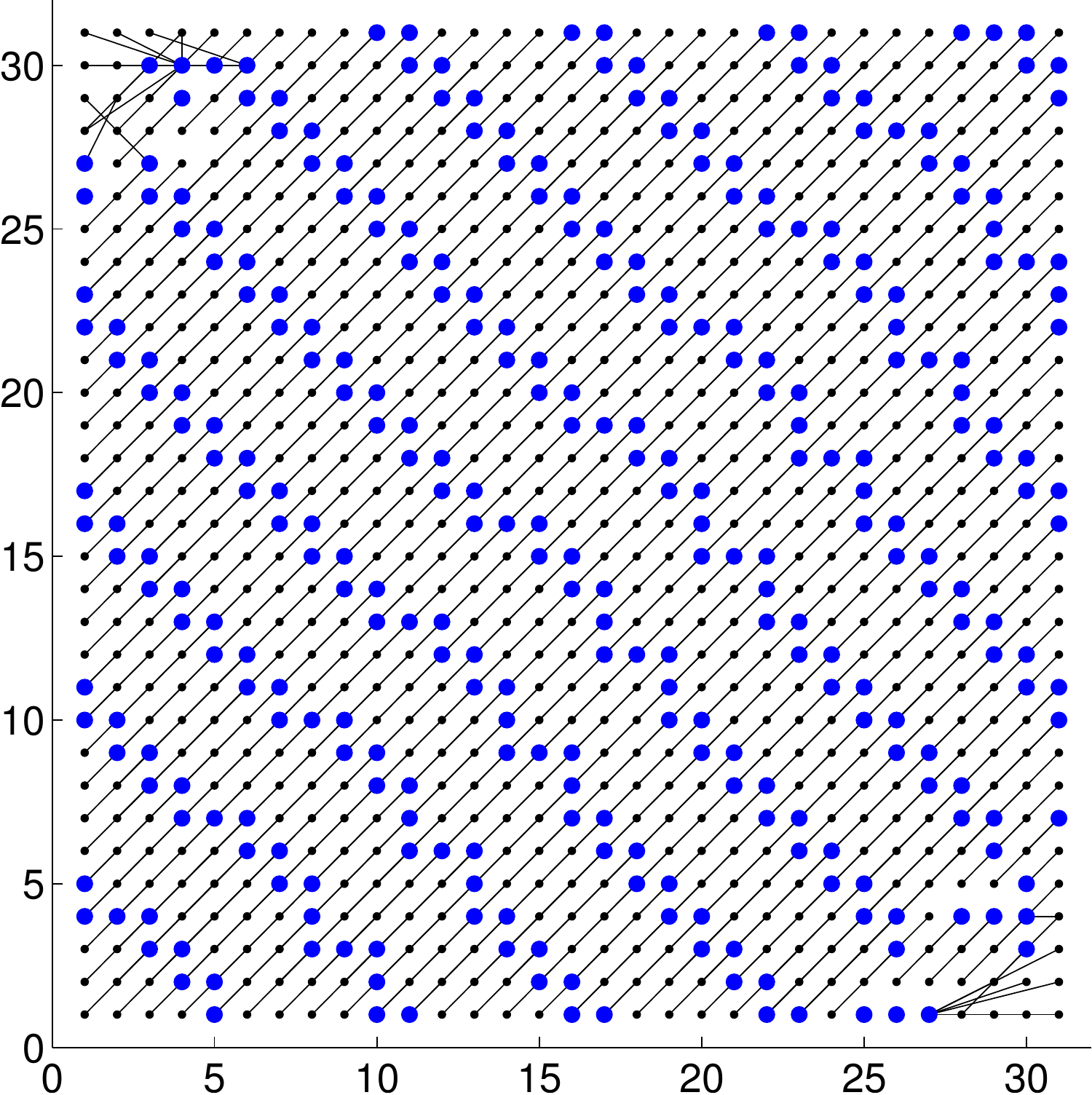}}
\hfill
  \subfigure[$\alpha = -\pi/4, \rho = .28, \rho_f = .45, \gamma_o = 1.862, \gamma_g = .469 $ \label{fig:s2mp4}]{\includegraphics[scale = 0.5]{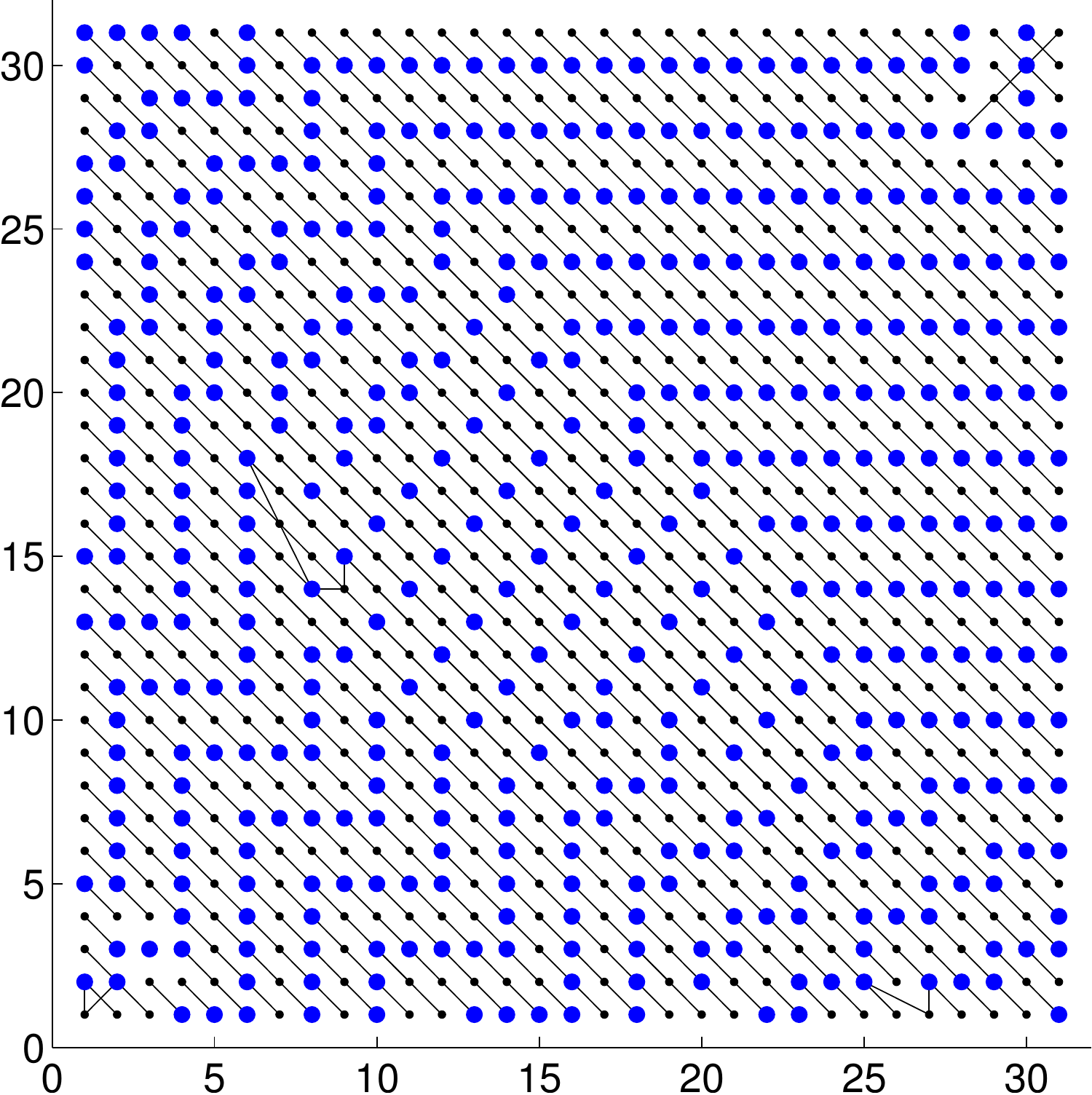}}
 \hfill
 \subfigure[$\alpha = \pi/8, \rho = .29, \rho_f = .59, \gamma_o = 1.805, \gamma_g = .295 $\label{fig:s2p16}]{\includegraphics[scale = 0.5]{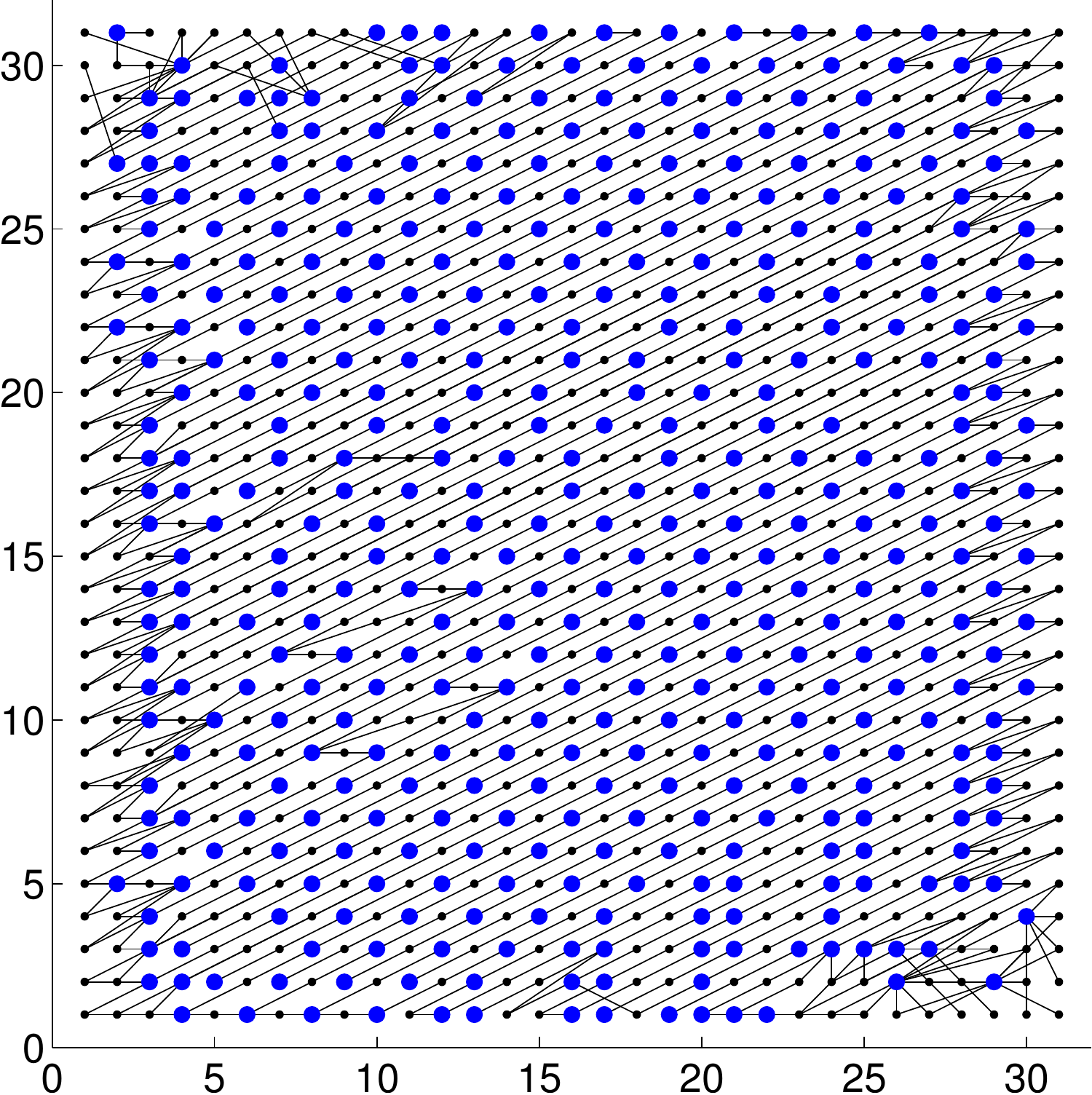}}  
    \caption{Coarse grids and interpolation patterns for $h= 1/32$ for various choices of $\alpha$, using the graph of $A^2$, i.e., $d  = 2$ and $d_{LS} = 4$,  to define the strength matrix. \label{fig:s2}}
\end{figure}

\begin{figure}
 \subfigure[$\alpha = \pi/8, \rho = .29, \rho_f = .59, \gamma_o = 1.805, \gamma_g = .295 $\label{fig:s2p6}]{\includegraphics[scale = 0.5]{figs/grid_api8_s2_rd4.pdf}}  \hfill
  \subfigure[$\alpha = \pi/8, \rho = .40, \rho_f = .76, \gamma_o = 1.543, \gamma_g = .257$\label{fig:s4p6}]{\includegraphics[scale = 0.5]{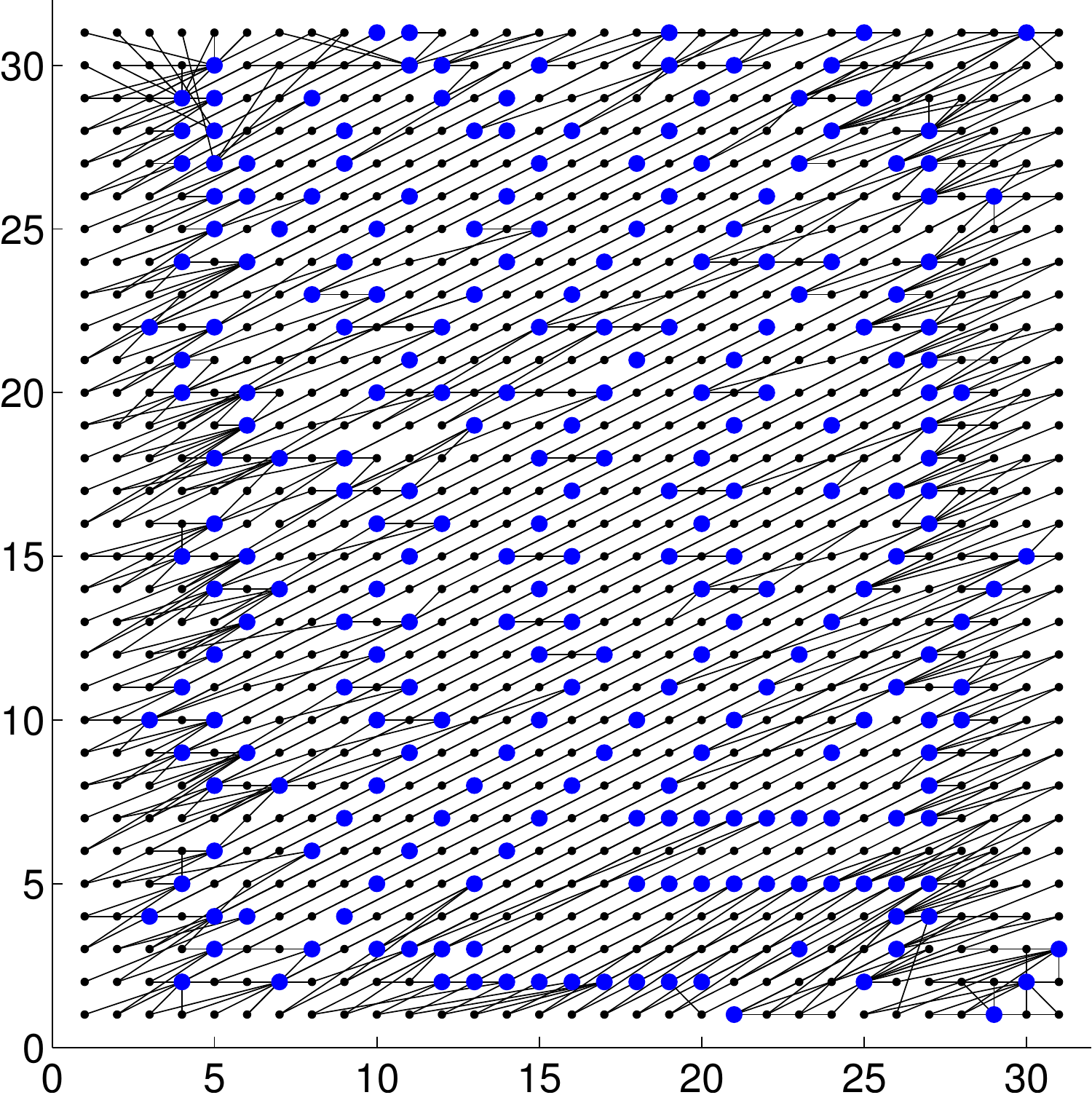}}
   \caption{Coarse grids and interpolation patterns for $h= 1/32$ with $\alpha=\pi/8$ using the graph of $A^2$, i.e., $d = 2$ and   $d_{LS}=4$ (left) and $A^4$, i.e., $d=4$ and $d_{LS} =6$ (right). 
\label{fig:s2k6}}
\end{figure}

In this section, we present the choices of coarse grids our algorithm makes when applied to the anisotropic problem for $\epsilon = 10^{-10}$ and various choices of the anisotropic angle $\alpha$ in \eqref{eq:a}. 
In the plots in Figures~\ref{fig:s1} -- \ref{fig:s2k6}, the black lines depict the interpolation pattern for each $\finevar$-point (denoted by the smaller  circles) from its neighboring $\coarsevar$-points (denoted by the larger  circles).  We note that for the $d=1$ tests
in Figure  \ref{fig:s1}, the coarsening and interpolation pattern  follow the anisotropy when the choice of the coarse grids allow it, i.e., for $\alpha = 0$ and $\alpha = \pi/4$. This is not the case for the $\alpha = -\pi/4$ and $\alpha =  \pi/8$ cases, for which the direction of strongest coupling is not captured by the connectivity of the fine-level discretization of the problem.   Hence, for these problems it is not possible to form a strength matrix from immediate algebraic neighbors which 
produces a coarse grid that allows  the interpolation pattern to 
follow the anisotropy.
These observations in fact led to the idea
of using the graph $\mathcal{G}_d$ of $A^d$, $d>1$, to form the strength
matrix.  Overall, we see that even for these cases, the algebraic
distance strength measure leads to least squares interpolation that follows the general direction of anisotropy as much as possible when coarsening is done using $\mathcal{G}_1$ to form the strength graph and set of coarse grid points.  

Another interesting observation here is that for the $d=2$ results in Figure \ref{fig:s2}, by using 
the graph of $A^2$ to form the algebraic-distance-based strength graph, the coloring algorithm is now able to coarsen in the direction of anisotropy for the $\alpha = -\pi/4$ case.  This choice also allows a more aggressive coarsening for other values of $\alpha$ while maintaining, for $\alpha = 0$ and $\alpha = \pi/4$, the characteristic directions induced by the anisotropy.

A particularly difficult choice of $\alpha$ occurs when the anisotropy direction is nearly aligned with a grid direction, for example, taking $\alpha \approx 0$ for our diffusion problems on a regular 2D grid.  For these cases, a longer-range interpolation and  an extended search depth for coarse-grid candidates is needed to accurately capture the anisotropy  (more so for directions closer to the grid lines).   To demonstrate this phenomenon, 
we include the results for $\alpha = \pi/8$ in Figure \ref{fig:s2k6}.  Here, we see that when using a larger value of $d$ in \eqref{eq:strength_graph} to define $\mathcal{G}_d$ and also a larger value of the search depth ($d_{LS}$ in \eqref{eq:nbhd}), the LS procedure is able to generate an interpolation operator that accurately follows the direction of the anisotropy.  


\subsection{The coarse-level operator}
One of the interesting  deliverables of the algorithm, in particular of its implementation of the compatible relaxation and the algebraic distances, is the pattern of the resulting coarse-grid stencil. Discretization involving 
(\ref{fd_seven})  favors $\alpha = \pi/4$ and with the same argument makes the worst possible discretization for  $\alpha = -\pi/4$ for which employment of  upper-left and lower-right  grid-point neighbors in discretization of $\partial_{xy}$ would be beneficial.     
We now consider both cases $\alpha = \pm \pi/4$ and  the seven point discretization, employing (\ref{fd_seven}). For both cases, we assume $\epsilon = 10^{-10}$,  $d = 2$, and $d_{LS} = 4$.

\begin{figure}
  \subfigure[Coarse-grid equation pattern for $\alpha = \pi/4$.\label{fig:good}]{\includegraphics[width=.4\textwidth]{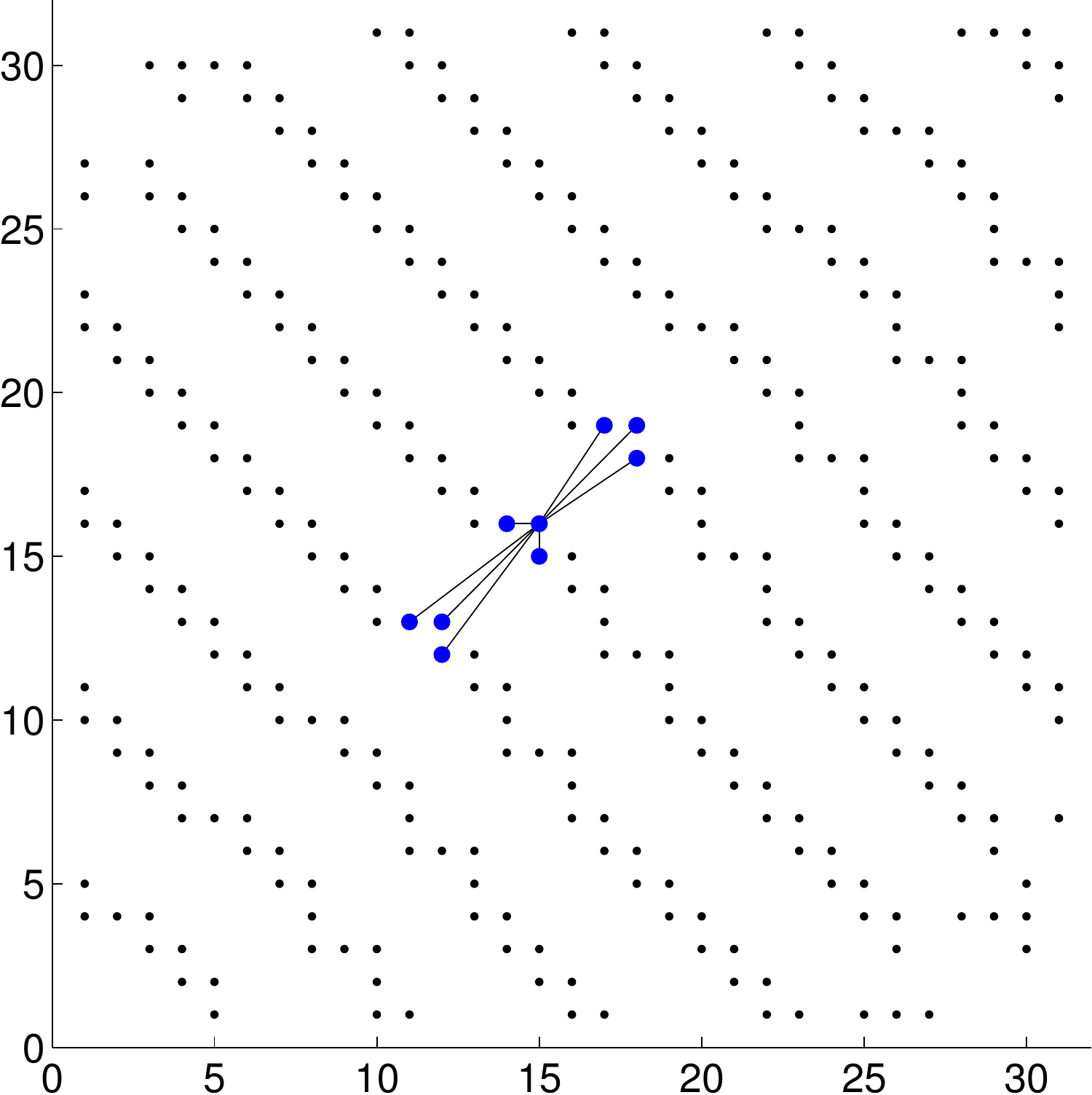}}\hfill
  \subfigure[Coarse-grid equation pattern for $\alpha = -\pi/4$ \label{fig:bad}]{\includegraphics[width=.4\textwidth]{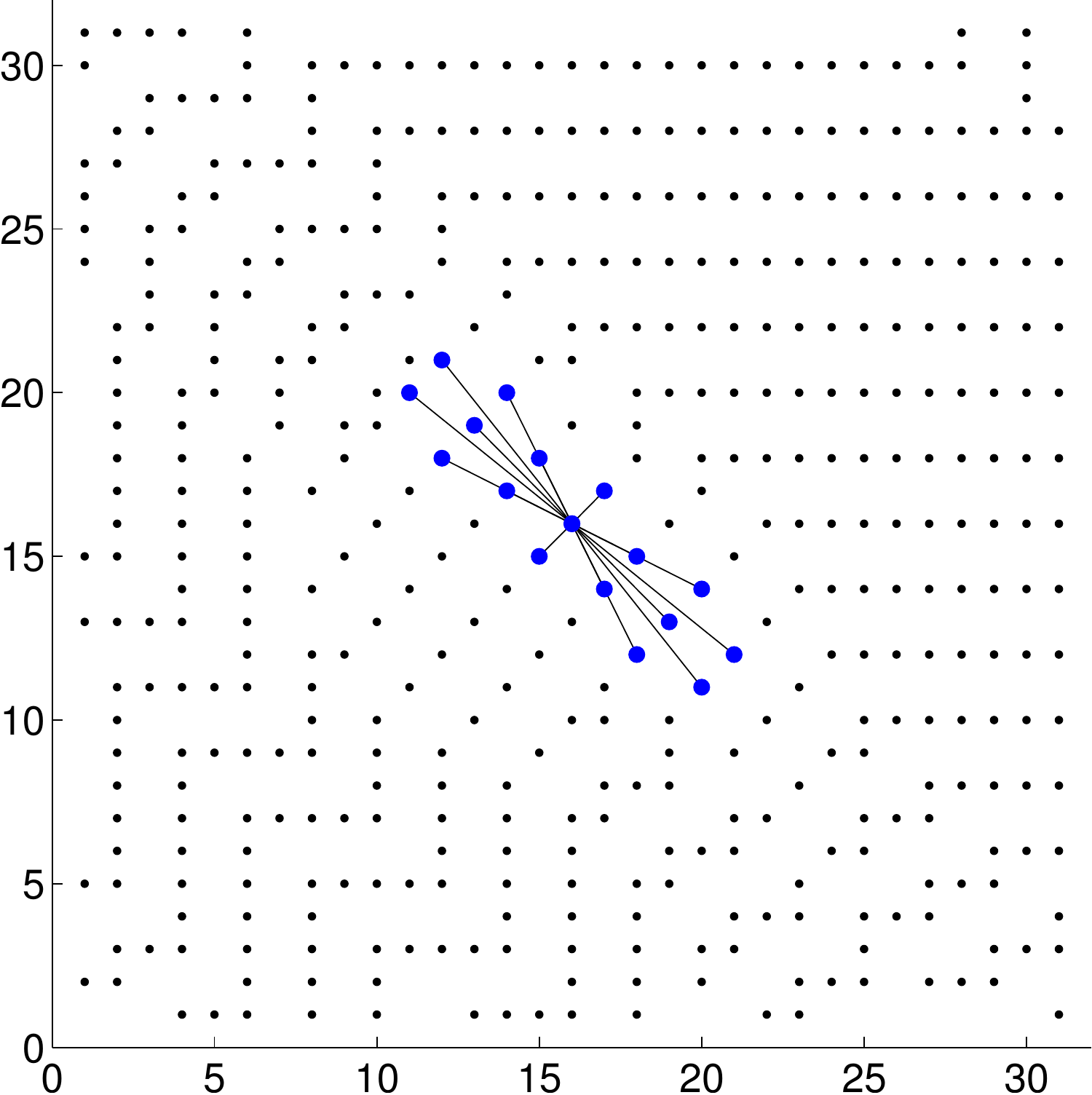}}
  \caption{ Non-zero pattern of one of the  stencils of the coarse-grid operator   $A_c$, centered at $i \in C$ and connected with 
  $j \in C$ such that $(A_c)_{ij} \neq 0$. The smaller dots in the graph are all other coarse-grid variables.    
  \label{fig:pis}}
\end{figure}

 We first confirm that for $\alpha = \pi/4$, the coarse-grid operator $A_c = P^t A P$ preserves the intrinsic strength of connections inherited from the fine-grid operator $A$.  A typical example of its stencil is  given in Fig \ref{fig:good} (the details of configurations depend on the coarsening pattern in the neighborhood of the considered coarse-grid equation).
 
 The results for the more challenging $\alpha = -\pi/4$ case are provided in Figure \ref{fig:bad}.  Here, we observe that although the discretization on the fine grid does not follow the anisotropy whatsoever, the non-zero pattern of the coarse-grid operator correctly aligns with the direction of anisotropy.  This result demonstrates the ability of the algorithm to overcome, if needed, the disadvantage of a poorly chosen fine-grid discretization and regain a more favorable discretization on the first coarse grid.  Further, the results for the $\alpha =  \pi/4$  indicate, that all consecutive coarse grids (though not employed in our two-level algorithm) are likely  to maintain a similar favorable discretization, that too accurately reflects the anisotropy.

Coefficients of the coarse-grid stencils,  presented in
Fig.\ref{fig:pis}, are given next. Here  $S_{cg}^{+}$ corresponds to
$\alpha = \pi/4$ (entries denoted by $*$  are  negligible,  with absolute values below $10^{-11}$) 
\[
S_{cg}^{+} = 
\begin{pmatrix}
     &   & & &  * &  -0.166\\
             &        & & &  & * &\\
           &  & * & 0.332 & &    \\
     &  &  & * & &  \\
        * & -0.166 &  &  & & \\
             & * & &  &  & 
                     \end{pmatrix}, 
    \]
    and $S_{cg}^{-} $ corresponds to  $\alpha = -\pi/4$:
    \small{     \[
                     S_{cg}^{-} = 
\left(
\begin{array}{*{12}{r}} 
            & -0.11 & & & & & & & &  &  &  \\
            -0.12 &  &  &  \phantom{-}0.23 & & & &  &   &  \\
            & & \phantom{-}1.37 & & & & & & &  &   &  \\
            &  0.23 & &  & -2.94  & & & &  &   & \\
           &  & &   -2.91 &   & & 1.06 & &   &  \\
            &  & &    &   &    \phantom{-}6.36 &   & &  &   &  \\
            &  & &    &  1.06  &  &   & -2.90 &        &   \\              
            &  & &    &  &  &  -2.93 &   &  &    0.23   &   \\           
                &  & &    &  &  &   &   &   \phantom{-} 1.38  &     &  \\
                       &  & &    &  &  &   &    0.23 &  & &  -0.12    \\ 
           &  & &    &  &  &   &    &  &  -0.12  &    
 \end{array} \right).
\]}
In both stencils, all  entries are rounded to the nearest hundredth.  

The distribution of the stencils' coefficients further illustrates the ability of the algebraic distance based 
strength of connections measure to choose the 
correct coarse-grid points for problems with anisotropic coefficients; in both cases, the non-zero sparsity pattern and dominant coefficients
of the resulting coarse-grid operators follow the direction of anisotropy.   

\begin{remark}
We recall that a larger graph distance $d =2$ and search depth $d_{LS} =4$ were required to obtain our promising results for the $\alpha = \pi/4$ case and note 
 the additional fill-in of the resulting coarse-level operator in this case, as seen in Figure~\ref{fig:pis}.b.  Generally, as the direction of the non-grid aligned anisotropy 
aligns itself more closely with the grid, the values of $d$ and, hence, $d_{LS}$, must be increased for the approach to appropriately capture the anisotropy.  This, in turn, 
increases the fill-in of the coarse-level operator, making it difficult to recursively coarsen the equations in a systematic way and maintain low grid and operator complexities.   We mention that any method for constructing the long-range interpolation required by such problems will have to be designed with this issue in mind.  
 \end{remark}

\subsection{Two-level convergence}
We conclude our experiments with tests of the proposed AMG setup algorithm applied to \eqref{eq:diff} for various choices of the anisotropy angle, $\alpha$, the anisotropy coefficient, 
$\epsilon$, and the mesh spacing $h = 1/N$.  In the following tables, the asymptotic convergence rates, $\rho$, of the two-grid solver produced by our setup algorithm  are reported, along with the corresponding coarsening factors $\gamma_g$ and  operator complexity ratios
 $\gamma_o$.  Here, we observe a slight dependence of the computed convergence rates and grid and operator complexities on the problem parameters, $\epsilon$, 
 $\alpha$, and $h$.  This dependence on $h$ is restricted mostly to the non-grid aligned cases, with the exception $\alpha = 0$ and $\epsilon = .1$, where we see a slight increase in $\rho$ as the problem size is increased from $N=64$ to $N=128$.   Moreover, in all cases, the convergence rates and complexities are uniformly bounded with respect to $\epsilon$ and $\alpha$ for fixed $h$.  We note in addition that these results are promising when considering that all tests were performed with the same strength parameter $\theta_{ad} = .5$.  In fact, all  parameters in the setup algorithm were fixed, suggesting that the individual components of the setup are robust for the targeted anisotropic problems, even those leading to non $M$-matrix systems as in the $\alpha = -\pi/4, \pi/8$ cases.  Further, we note that the setup handles the isotropic case when $\alpha = 0$ and $\epsilon =1$ with similar efficiency, producing a two-grid method with convergence rate $\rho = .28$ and complexities $\gamma_g = 25$ and $\gamma_o = 1.6$ for $N=32,64,128$.

\begin{table}[!ht]
  \begin{center}
  
  \begin{tabular}{|c|c|c|c|c|c|c|c|}
      \hline
      $\alpha$\\
      \hline
      \hline
       $\pi/4$\\
      \hline
       $-\pi/4$\\
       \hline
       $\pi/8$ \\ 
       \hline
       $0$  \\
        \hline
         \end{tabular}
    \begin{tabular}{|c|c|c|c|c|c|c|c|}
    \hline
        $\epsilon =  .1$ \\ 
        \hline      \hline
         .10 (.35,1.5) / .22 (.30,1.5)  / .22 (.27,1.5) \\
       \hline 
          .31 (.35,1.7) / .36 (.34,1.7)  / .48 (.32,1.7)   \\
       \hline
          .32 (.27,1.4) /  .39 (.24,1.5)  / .35 (.25,1.5)  \\
        \hline
         .19 (.34,1.6) /  .20 (.37,1.7) / .24 (.38,1.7)  \\
      \hline
         \end{tabular}
\bigskip

 \begin{tabular}{|c|c|c|c|c|c|c|c|}
      \hline
      $\alpha$\\
      \hline
      \hline
       $\pi/4$\\
      \hline
       $-\pi/4$\\
       \hline
       $\pi/8$ \\ 
       \hline
       $0$  \\
        \hline
         \end{tabular}
    \begin{tabular}{|c|c|c|c|c|c|c|c|}
    \hline
    $\epsilon = .0001$ \\
        \hline      \hline
         .26 (.35,1.5) /  .26 (.34,1.5)  / .23 (.36,1.5)  \\
       \hline 
        .28 (.46,1.9) /   .33 (.45,1.9)  /  .38 (.41,1.9)   \\
       \hline
           .30 (.43,1.8) / .48 (.38,1.8) / .51 (.36,1.8)    \\
        \hline
        .05 (.34,1.6) /  .06 (.35,1.6)   / .06 (.38,1.7)   \\
      \hline
         \end{tabular} 
        \bigskip

 \begin{tabular}{|c|c|c|c|c|c|c|c|}
      \hline
      $\alpha$\\
      \hline
      \hline
       $\pi/4$\\
      \hline
       $-\pi/4$\\
       \hline
       $\pi/8$ \\ 
       \hline
       $0$  \\
        \hline
         \end{tabular}
    \begin{tabular}{|c|c|c|c|c|c|c|c|}
    \hline
       $\epsilon = 0$\\
        \hline      \hline
         .06 (.33,1.4) / .06 (.34,1.4)  / .06 (.36,1.5)\\
       \hline 
       .28 (.46,1.9) / .35 (.43,1.9) / .37 (.41,1.9)    \\
       \hline
         .30 (.43,1.8) / .49 (.38,1.8) / .52 (.36,1.8) \\
        \hline
         .05 (.34,1.3) /  .06 (.35,1.3)  / .06 (.38,1.4)  \\
      \hline
         \end{tabular}
         \bigskip

   \caption[Two-grid LS interpolation for FD Laplacian]{Approximate asymptotic convergence rates for the seven point FD Anisotropic Laplace problem with Dirichlet boundary conditions for various choices of $\alpha$, $\epsilon$ and $h$.  Here, our proposed setup is applied with $d=2$ and $d_{LS} = 4$.  The reported results correspond to: $\rho$ ($\gamma_g,\gamma_0)$ for $n = 32^2$ /  $n=64^2$ / $n=128^2$.}  
  \end{center}
\end{table}

\section{Concluding remarks}
\label{sec:conclusions}
The LS functional gives a flexible and robust tool for measuring AMG strength of connectivity via algebraic distances:
between pairs of points to define a strength graph used to choose coarse points; and
among sets of points for determining interpolatory sets.  
The proposed coarsening approach combining algebraic distances, compatible relaxation, and least squares interpolation provides an effective scheme for the non-grid aligned anisotropic diffusion problems considered. 
The approach chooses suitable coarse-grid variables and prolongation operators for a wide range of anisotropies, without the need for parameter tuning.  
 Moreover, even when the initial fine-grid discretization is chosen in the direction opposite to the one defined by the  anisotropy (as in the $\alpha = -\pi/4$ case), the method constructs a suitable interpolation operator and, further, produces a coarse-grid operator which better captures the anistropy. While not the focus of this paper, the promising results obtained by our approach suggest that its extension from two to many levels should be effective. Indeed, it is natural to assume that if the first coarse grid system produced by our scheme is consistent with a suitable discretization, regardless of the initial discretization, then the multilevel scheme constructed from it will yield an effective solver.  As noted earlier, the main challenge faced in extending our approach to a multilevel one for the targeted anisotropic problems is that of designing an algorithm capable of constructing long-range interpolation as needed to capture general anisotropies and at the same time maintains low grid and operator complexities.
\bibliographystyle{alpha}
\bibliography{bamg_Aniso}
\end{document}